\documentclass{ensmath}
\title[CONNECTIVITY OF COMPLEXES ASSOCIATED TO SURFACES]{A NOTE ON THE CONNECTIVITY OF CERTAIN COMPLEXES ASSOCIATED TO SURFACES}
\author[A. PUTMAN]{Andrew \sn{Putman}}
\usepackage{amsmath}
\usepackage{amssymb}
\usepackage{epsfig}
\usepackage{amsfonts}
\usepackage{amscd}

\newcommand\Heading[1]{{\bf {\smallskip}{\noindent}#1 :}}

\newcommand\InsertFigure[3]{
\begin{figure}[tbp]
\centering
\includegraphics[scale=0.75]{#2}
\caption{#3}
\label{#1}
\end{figure}}

\newcommand\CaptionSpace{\hspace{0.2in}}

\newcommand\Mod{\mbox{Mod}}
\newcommand\Torelli{\mbox{${\mathcal I}$}}
\newcommand\Z{\mbox{$\mathbb{Z}$}}
\newcommand\HHHH{\mbox{H}}
\newcommand\Curves{\mbox{$\mathcal{C}$}}
\newcommand\CSep{\mbox{$\mathcal{C}_{\text{sep}}$}}
\newcommand\CNosep{\mbox{$\mathcal{C}_{\text{nosep}}$}}
\newcommand\CHalf{\mbox{$\mathcal{C}_{\text{half}}$}}
\newcommand\CutSystem{\mbox{$\mathcal{CT}$}}
\newcommand\Pants{\mbox{$\mathcal{P}$}}

\newcommand\CSepk[1]{\mbox{$\mathcal{C}_{\text{sep}}^{(#1)}$}}

\newcommand\CHalfk[1]{\mbox{$\mathcal{C}_{\text{half}}^{(#1)}$}}
\newcommand\CutSystemk[1]{\mbox{$\mathcal{CT}^{(#1)}$}}
\newcommand\Pantsk[1]{\mbox{$\mathcal{P}^{(#1)}$}}
\newcommand\Cayley{\mbox{Cay}}

\begin{document}

\maketitle

\begin{abstract}
This note is devoted to a trick which yields almost trivial
proofs that certain complexes associated to topological surfaces are connected
or simply connected.  Applications include new proofs that the
complexes of curves, separating curves, nonseparating curves, pants,
and cut systems are all connected for genus $g \gg 0$.  We also prove that
two new complexes are connected : one involves curves which split
a genus $2g$ surface into two genus $g$ pieces, and the other
involves curves which are homologous to a fixed curve.  The connectivity
of the latter complex can be interpreted as saying the ``homology''
relation on the surface is (for $g \geq 3$) generated
by ``embedded/disjoint homologies''.  We finally prove that the
complex of separating curves is simply connected for $g \geq 4$.
\end{abstract}

\section{Introduction}
Let $\Sigma_{g}$ be a genus $g$ surface and 
$\Mod(\Sigma_{g})$ be the mapping class group of $\Sigma_{g}$, that is, the
group of isotopy classes of orientation-preserving homeomorphisms of $\Sigma_{g}$ (see
\cite{IvanovSurvey} for a survey of $\Mod(\Sigma_g)$).
An important theme in the study of $\Mod(\Sigma_g)$ and its subgroups is the close
relationship between algebraic properties of $\Mod(\Sigma_g)$ (e.g.\ cohomology, finiteness
properties, automorphisms, etc.)\ and the structure of $1$-submanifolds of $\Sigma_g$.  The combinatorics
of these $1$-submanifolds have been encoded in the structure of a number
of simplicial complexes, such as the curve complex and the pants complex.  A key property
of these complexes is that they are often highly connected.  In this paper, we
discuss a general trick which yields simple proofs that complexes of this sort
are connected or simply connected; in many cases this is sufficient for the applications.

In the past, these sorts of theorems have been proven using a variety of tools, such as
curve surgery (see, e.g., \cite{FarbIvanov, HatcherEasy, MasurSchleimer, McCarthyVautaw, SchleimerNotes, WajnrybElementary}), 
parametrized Morse theory (see, e.g., \cite{HatcherThurston, IvanovCurves}), and 
Teichm\"{u}ller theory (see, e.g., \cite{BowditchEpstein, HarerVirtual, PennerTeichmuller}).  We
instead exploit the basic combinatorial group theoretic properties of $\Mod(\Sigma_{g})$
and its subgroups, deducing that complexes are connected from the structure of generating sets 
and deducing that they are simply connected from relations.  Of course, we may be
accused of circular reasoning, as the standard construction of generators and relations
for $\Mod(\Sigma_g)$ involves investigating connected and simply connected complexes upon which it acts!
The point of this paper is that this only needs to be done once -- as soon as generators and relations for
$\Mod(\Sigma_g)$ are found, one can prove that essentially any complex upon which $\Mod(\Sigma_g)$ acts in a reasonable way is connected or simply
connected by a formal, finitely checkable (and in practice quite easy) process.

The first complex we will examine is the complex of curves (introduced by Harvey in \cite{HarveyCurves}), together
with two of its subcomplexes.

\begin{definition}
The {\em complex of curves} $\Curves(\Sigma_g)$ is the simplicial complex whose
simplices are sets $\{c_1,\ldots,c_k\}$ of non-trivial isotopy classes of simple closed
curves on $\Sigma_g$ which can be realized disjointly.  The {\em complex of separating curves} $\CSep(\Sigma_g)$ and the 
{\em complex of nonseparating curves} $\CNosep(\Sigma_g)$ 
are the full subcomplexes of $\Curves(\Sigma_g)$ spanned by separating and nonseparating curves, respectively. 
\end{definition}

\noindent
We will give a simple, unified proof of the following theorem, which for $\Curves(\Sigma_g)$ and
$\CNosep(\Sigma_g)$ is due to Lickorish \cite{LickorishGenerators} (though he did not use this language) and
for $\CSep(\Sigma_g)$ is due to Farb and Ivanov \cite{FarbIvanov}.  Other proofs
of the connectedness of $\CSep(\Sigma_g)$ can be found in \cite{MasurSchleimer} and \cite{McCarthyVautaw}.

\begin{theorem}
\label{theorem:curvesconnected}
$\Curves(\Sigma_g)$ and $\CNosep(\Sigma_g)$ are connected for $g \geq 2$, while $\CSep(\Sigma_g)$
is connected for $g \geq 3$.
\end{theorem}

In fact, our trick allows us to achieve rather precise control over the topology 
of the curves which appear in our complexes.  For instance, consider the following complex.

\begin{definition}
Let $\CHalf(\Sigma_{2g})$ be the simplicial complex whose simplices are sets $\{c_1,\ldots,c_k\}$ of
isotopy classes of simple closed curves on $\Sigma_{2g}$ which satisfy the following two conditions.
\begin{itemize}
\item Each $c_i$ separates $\Sigma_{2g}$ into two genus $g$ subsurfaces.
\item For $i \neq j$, the geometric intersection number $i_g(c_i,c_j)$ is {\it minimal} among such curves.  This minimality
means the following : if $g = 1$, then $i_g(c_i,c_j)=4$, while if $g \geq 2$, then $i_g(c_i,c_j)=2$.
\end{itemize}
\end{definition}

\noindent
We will prove the following theorem, answering a question posed to the author by Schleimer (who
proved the theorem for $g = 1$ \cite{SchleimerNotes}).

\begin{theorem}
\label{theorem:halfconnected}
$\CHalf(\Sigma_{2g})$ is connected for $g \geq 1$.
\end{theorem}

We next investigate the cut system and pants graphs, which were introduced
by Hatcher and Thurston in \cite{HatcherThurston}.  

\begin{definition}
A {\em cut system} on $\Sigma_g$ is a set $\{c_1,\ldots,c_g\}$ of isotopy classes of simple closed
curves on $\Sigma_g$ which can be realized disjointly with $\Sigma_g \setminus (c_1 \cup \cdots \cup c_g)$ connected (see
Figure \ref{figure:modconnected}.d).  Two cut systems $\{c_1,\ldots,c_g\}$ and $\{c_1',\ldots,c_g'\}$ 
differ by an {\em elementary move} if
there is some $1 \leq i \leq k$ so that $i_g(c_i,c_i')=1$ and so that $c_j = c_j'$ for $j \neq i$.  The
{\em cut system graph} $\CutSystem(\Sigma_g)$ is the graph whose vertices are cut systems on $\Sigma_g$
and whose edges correspond to elementary moves between cut systems.
\end{definition}

\begin{definition}
For $g \geq 2$, a {\em pants decomposition} of $\Sigma_g$ is a maximal simplex $\{c_1,\ldots,c_k\}$ of $\Curves(\Sigma_g)$ (see
Figure \ref{figure:modconnected}.e).  Observe that $k = 3g-3$ and that cutting $\Sigma_g$ along the 
$c_i$ results in a collection of $3$-holed spheres (the ``pairs of pants'').
Two pants decompositions $\{c_1,\ldots,c_k\}$ and $\{c_1',\ldots,c_k'\}$ differ by an {\em elementary move} if
there is some $1 \leq i \leq k$ so that for $j \neq i$ we have $c_j = c_j'$ and so that $i_g(c_i,c_i')$ is {\it minimal} among such
curves.  This minimality means the following (see Figure \ref{figure:modconnected}.g) : if $S$ is the component of $\Sigma_g$ cut along 
$c_1 \cup \cdots \cup c_{i-1} \cup c_{i+1} \cup \cdots \cup c_k$ containing $c_i$, then $i_g(c_i,c_i')=2$ if $S$ is a 4-holed sphere and
$i_g(c_i,c_i')=1$ if $S$ is a 1-holed torus.  The {\em pants graph} $\Pants(\Sigma_g)$ is the graph whose vertices
are pants decompositions of $\Sigma_g$ and whose edges correspond to elementary moves between pants decompositions.
\end{definition}

\begin{remark}
Hatcher and Thurston in fact considered $\CutSystem(\Sigma_g)$ and $\Pants(\Sigma_g)$ with a number of 2-cells
attached to render them simply connected.  We will make no use of these 2-cells.
\end{remark}

\noindent
We will give a new proof of the following theorem of Hatcher and Thurston, which for $\CutSystem(\Sigma_g)$ is
Theorem 1.1 of \cite{HatcherThurston} and for $\Pants(\Sigma_g)$ is contained in the appendix of \cite{HatcherThurston}.

\begin{theorem}
\label{theorem:cutpantsconnected}
$\CutSystem(\Sigma_g)$ is connected for $g \geq 1$, while $\Pants(\Sigma_g)$ is connected for $g \geq 2$.
\end{theorem}

Next, we will use the action of the Torelli subgroup of the mapping class group (defined below)
to prove the following theorem, which elucidates the nature of the homology relation
on a surface.  It says that this relation is generated by ``embedded homologies'' (in the statement
of this theorem and throughout this paper, when we say that two simple closed unoriented curves
are homologous, we mean that they can be oriented in such a way that they are rendered homologous).

\begin{theorem}
\label{theorem:homologyconnected}
Fix $g \geq 3$, and let $\gamma$ and $\gamma'$ be homologous non-trivial simple closed curves on
$\Sigma_g$.  Then there exists a sequence
$$\gamma = \gamma_1, \gamma_2, \ldots, \gamma_k = \gamma'$$
of non-trivial simple closed curves on $\Sigma_g$ so that for $1 \leq i < k$ the curves $\gamma_i$
and $\gamma_{i+1}$ are disjoint and there exists an embedded subsurface $S_i \hookrightarrow \Sigma_g$
with $\partial S_i = \gamma_i \sqcup \gamma_{i+1}$ (in particular, $\gamma_i$ and $\gamma_{i+1}$ are homologous).
\end{theorem}

\begin{remark}
This theorem is false for $g=2$, as there exist no subsurfaces $S$ of $\Sigma_2$ so that $\partial S$ consists
of two simple closed curves which are nonseparating and nonisotopic on $\Sigma_2$.
\end{remark}

Finally, we will show that our methods can be extended to prove that various complexes are simply connected.  As
an example, we prove the following.

\begin{theorem}
\label{theorem:csepsimplyconnected}
For $g \geq 4$, the complex $\CSep(\Sigma_g)$ is simply connected.
\end{theorem}

\begin{remark}
Hatcher and Vogtmann \cite{HatcherVogtmann} have proven a much stronger theorem which says 
that $\CSep(\Sigma_g)$ is $\lfloor \frac{g-3}{2} \rfloor$-connected.  Their
result, however, does not imply Theorem \ref{theorem:csepsimplyconnected} for $g=4$.
\end{remark}

\begin{notation}
Let $P_1,P_2\ldots,P_k$ be a sequence of paths in a simplicial complex $X$ each of which begins and ends 
in the $0$-skeleton $X^{(0)}$ (we allow degenerate paths $P_i$ consisting of single vertices).  For all
$1 \leq i < k$, let $q_i$ be the terminal point of $P_i$ and $p_{i+1}$ be the initial point of $P_{i+1}$, and assume
that $\{q_i,p_{i+1}\} \in X^{(1)}$.  Thus either $q_i = p_{i+1}$ or $\{q_i,p_{i+1}\}$ is a 1-simplex.  We then denote
the path which first traverses $P_1$, then $P_2$, etc.\ by
$$P_1 - P_2 - \ldots - P_k.$$
\end{notation}

\section{Connectivity}
Our trick for proving that complexes are connected is contained in the following easy lemma.
\begin{lemma}
\label{lemma:main}
Consider a group $G$ acting upon a simplicial complex $X$.  Fix
a basepoint $v \in X^{(0)}$ and a set $S$ of generators
for $G$.  Assume the following hold.
\begin{enumerate}
\item For all $v' \in X^{(0)}$, the orbit 
$G \cdot v$ intersects the connected component of $X$ containing $v'$. 
\item For all $s \in S^{\pm 1}$, there is some path
$P_s$ in $X$ from $v$ to $s \cdot v$.
\end{enumerate}
Then $X$ is connected.
\end{lemma}
\begin{proof}
Consider $v' \in X^{(0)}$.  By Condition 1, there is some $g \in G$
together with a path $P$ from $g \cdot v$ to $v'$.  Write
$g$ as a word $s_1 \cdots s_k$ in $S^{\pm 1}$.  Then
$$P_{s_1} - s_1 P_{s_2} - \ldots - s_1 s_2 \cdots s_{k-1} P_{s_k} - P$$
is a path from $v$ to $v'$. \qed
\end{proof}

\InsertFigure{figure:modconnected}{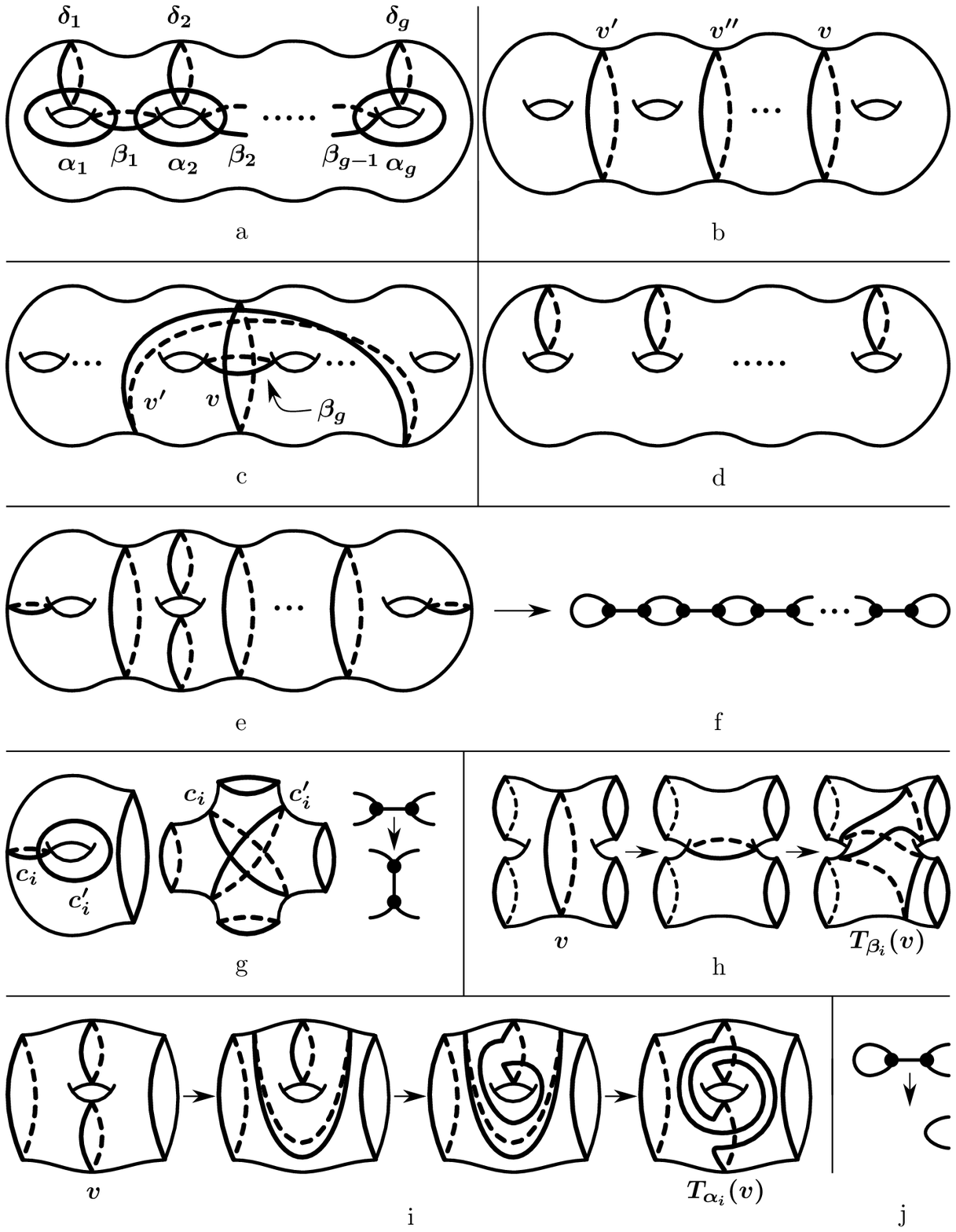}{The various figures needed for the proof of Theorems \ref{theorem:curvesconnected}, 
\ref{theorem:halfconnected}, and \ref{theorem:cutpantsconnected} (see that proof for more details)}

\noindent
We now prove Theorems \ref{theorem:curvesconnected}, \ref{theorem:halfconnected}, and \ref{theorem:cutpantsconnected}.

\begin{likeproof}[Proof of Theorems \ref{theorem:curvesconnected}, \ref{theorem:halfconnected}, and \ref{theorem:cutpantsconnected}]
Let 
$$S=\{T_{\alpha_1},T_{\delta_1},\ldots,T_{\alpha_g},T_{\delta_g},T_{\beta_1},\ldots,T_{\beta_{g-1}}\}$$ 
be the Dehn twists about the curves in Figure \ref{figure:modconnected}.a; Lickorish proved that $S$ generates 
$\Mod(\Sigma_g)$ (see \cite[\S 4]{IvanovSurvey} for the definition of a Dehn twist and a discussion of Lickorish's 
theorem).  For each complex in question, we will pick a basepoint $v$ and verify that the two conditions of Lemma 
\ref{lemma:main} are satisfied for the action of the mapping class group upon the complex.  We emphasize
that in each proof the basepoint $v$ and all other ancillary curves are chosen to intersect the (already fixed)
generators in simple ways.  In particular, they depend on the choice of generators.

\Heading{$\boldsymbol{\Curves(\Sigma_g)}$, $\boldsymbol{\CNosep(\Sigma_g)}$, and $\boldsymbol{\CSep(\Sigma_g)}$}
The proofs for these three complexes are similar; we will discuss $\CSep(\Sigma_g)$ and leave the other two to
the reader.  Our basepoint $v \in \CSepk{0}(\Sigma_g)$ will be the curve indicated in Figure \ref{figure:modconnected}.b.  
The orbit $\Mod(\Sigma_g) \cdot v$ consists of all separating curves which cut off 1-holed
tori.  Indeed, let $w$ be another separating curve which cuts off a 1-holed torus.  The classification of surfaces
implies that we get homeomorphic non-connected surfaces when we cut $\Sigma_g$ along either $v$ or $w$.  Gluing
together homeomorphisms between the cut surfaces yields the desired homeomorphism of 
$\Sigma_g$ taking $v$ to $w$ (this trick will be used repeatedly; we will call
it the {\em classification of surfaces trick}).  
Since every separating curve is adjacent (in $\CSep(\Sigma_g)$) to a curve which cuts off a 1-holed torus, condition 1 
follows.  To check condition 2, we will need the curve $v'$ indicated in Figure \ref{figure:modconnected}.b.
Consider $s \in S^{\pm 1}$.  If $s = T_{\beta_{g-1}}^{\pm 1}$, then $v - v' - s \cdot v$ is the desired path.  
Otherwise, we have $s \in S^{\pm 1}$ but $s \neq T_{\beta_{g-1}}^{\pm 1}$, so $s \cdot v = v$.  Condition 2 follows.

\Heading{$\boldsymbol{\CHalf(\Sigma_{2g})}$}
Here $S$ is the corresponding generating set for $\Mod(\Sigma_{2g})$.  
Our basepoint $v \in \CHalfk{0}(\Sigma_{2g})$ will be the curve indicated in Figure \ref{figure:modconnected}.c.  
If $g \geq 2$, we will also
need the ancillary curve $v'$ from the same figure.  Now, by the classification of surfaces trick, 
$\Mod(\Sigma_{2g})$ acts transitively on $\CHalfk{0}(\Sigma_{2g})$,
so condition 1 is trivial.  To check condition 2, consider $s \in S^{\pm 1}$.  If $s = T_{\beta_g}^{\pm 1}$, then for $g=1$ the vertices
$v$ and $s \cdot v$ are adjacent, while for $g \geq 2$, the desired path is $v - v' - s \cdot v$.  
If instead $s \neq T_{\beta_g}^{\pm 1}$, then $s \cdot v = v$. Condition 2 follows.

\Heading{$\boldsymbol{\CutSystem(\Sigma_g)}$}
Our basepoint $v \in \CutSystemk{0}(\Sigma_g)$ will be the cut system indicated in Figure \ref{figure:modconnected}.d.  By
the classification of surfaces trick, $\Mod(\Sigma_g)$ acts transitively on $\CutSystemk{0}(\Sigma_g)$, so condition 1 holds.  
Also, for $s \in S^{\pm 1}$, either $s \cdot v = v$ or $s \cdot v$ is adjacent to $v$, so condition 2 holds.

\Heading{$\boldsymbol{\Pants(\Sigma_g)}$}
Our basepoint $v \in \Pantsk{0}(\Sigma_g)$ will be the pants decomposition indicated in Figure \ref{figure:modconnected}.e.  We start by
verifying condition 2.  Consider $s \in S^{\pm 1}$.  If $s = T_{\delta_i}^{\pm 1}$, then $s \cdot v = v$.  If
$s = T_{\beta_i}$, then Figure \ref{figure:modconnected}.h contains the desired path (we only draw the portion
of the pants decomposition which changes).  A similar path works if $s = T_{\beta_i}^{-1}$.  
If $s = T_{\alpha_1}^{\pm 1}$ or $s = T_{\alpha_g}^{\pm 1}$, then
$s \cdot v$ is adjacent to $v$.  If $s = T_{\alpha_i}$ but $i \neq 1$ and $i \neq g$, then Figure \ref{figure:modconnected}.i contains
the desired path.  A similar path works if $s = T_{\alpha_i}^{-1}$ with $i \neq 1$ and $i \neq g$.  Condition 2 follows.

We now verify condition 1.  It is enough to show that $\Pants(\Sigma_g) / \Mod(\Sigma_g)$ is connected.  For each 
pants decomposition $p=\{c_1,\ldots,c_k\}$ of $\Sigma_g$, define a graph $\phi(p)$ as follows (see
Figure \ref{figure:modconnected}.f).  The vertices of $\phi(p)$ are the connected components of $\Sigma_g$ cut along 
the $c_i$ (the ``pairs of pants'').  The edges are in bijection with the curves $c_i$; the edge corresponding
to $c_i$ connects the vertices corresponding to the components on either side of $c_i$.  Thus $\phi(p)$ is
a trivalent graph with $2g-2$ vertices (a loop at a vertex counts as $2$ edges abutting that vertex). It
is clear that each such graph comes from a pants decomposition.  Moreover, it is not hard to see that
for pants decompositions $p$ and $p'$ we have $\phi(p)$ isomorphic to $\phi(p')$ if and only if there is some $f \in \Mod(\Sigma_g)$
so that $p = f \cdot p'$.

Now consider an elementary move from $p=\{c_1,\ldots,c_k\}$ to $p'=\{c_1',\ldots,c_k'\}$.
Let $c_i$ be the curve which changes in this move.  If $i_g(c_i,c_i')=1$ (so $c_i$ corresponds to a loop
in $\phi(p)$; see the left hand part of Figure \ref{figure:modconnected}.g), then $\phi(p)=\phi(p')$.  If $i_g(c_i,c_i')=2$ (see
the central part of Figure \ref{figure:modconnected}.g), then $\phi(p)$ is transformed into $\phi(p')$ in the following
way (see the right part of Figure \ref{figure:modconnected}.g) : we first collapse the edge in $\phi(p)$
corresponding to $c_i$, yielding a vertex of valence 4, which we then ``expand'' to two vertices of valence 3, each
of which abuts 2 of the edges which once abutted the vertex of valence 4.  We will call this an {\em elementary shift} of the graph.
It is not hard to see that any elementary shift of $\phi(p)$ is induced by an elementary move of $p$.

It is enough, therefore, to prove that if $G$ and $G'$ are trivalent graphs with the same (necessarily even) 
number of vertices, then $G$ may be transformed into $G'$ by a sequence of elementary shifts.  
The proof will be by induction on the number $k$ of vertices.  The base case $k=2$ being trivial, we
assume that $k > 2$.  Since neither $G$ nor $G'$ can be a tree, each must contain
a simple closed edge-path.  Transform $G$ and $G'$ by elementary shifts so that these closed edge paths
are as short as possible.  Observe that these minimal-length closed edge paths must be loops --
if they were not loops, then we could shorten them by performing elementary shifts which collapse edges in them.  
Let $\overline{G}$ and $\overline{G}'$ be the result of removing these loops, deleting the resulting
valence 1 vertices, and then finally deleting the resulting valence 2 vertices while combining the 2 edges abutting
them into a single edge (see Figure \ref{figure:modconnected}.j).  By induction we can convert $\overline{G}$ into $\overline{G}'$ by a sequence of elementary
shifts.  It is easy to see that we can then ``lift'' this sequence of elementary shifts to $G$, thus proving the
theorem. \qed
\end{likeproof}

Next, we prove Theorem \ref{theorem:homologyconnected}.

\begin{likeproof}[Proof of Theorem \ref{theorem:homologyconnected}]
This theorem is clearly equivalent to the connectedness of the following complex for $g \geq 3$.

\begin{definition}
Let $\Curves^{\gamma}(\Sigma_g)$ denote the full subcomplex of $\Curves(\Sigma_g)$ spanned
by curves homologous to $\gamma$.
\end{definition}

If $\gamma$ is separating, then $\Curves^{\gamma}(\Sigma_g) = \CSep(\Sigma_g)$, which is connected by 
Theorem \ref{theorem:curvesconnected}.  Assume, therefore, that
$\gamma$ is nonseparating, and let $\Torelli(\Sigma_g) \subset \Mod(\Sigma_g)$ (the {\em Torelli group}) be 
the kernel of the action of $\Mod(\Sigma_g)$ on $\HHHH_1(\Sigma_g;\Z)$.  We will apply Lemma \ref{lemma:main} to the
action of $\Torelli(\Sigma_g)$ on $\Curves^{\gamma}(\Sigma_g)$.  

To apply Lemma \ref{lemma:main}, we need a base point and a generating set.  Since $\gamma$ is nonseparating, the classification
of surfaces trick implies that there is a homeomorphism taking $\gamma$ to the curve $v$ depicted in Figure
\ref{figure:boundingtwist}.a.  We can therefore assume without loss of generality that $\gamma$ in fact equals
the curve $v$; this will be our base point.  It is well known (see, e.g., \cite[Lemma 6.2]{PutmanCutPaste})
that $\Torelli(\Sigma_g)$ acts transitively on the $0$-skeleton of $\Curves^{\gamma}(\Sigma_g)$, so condition 1 is trivial.

For the generating set, recall that Johnson \cite{JohnsonFinite} proved
that $\Torelli(\Sigma_g)$ is finitely generated.  Our generating set $S$ will be the generating set for $\Torelli(\Sigma_g)$ 
constructed in \cite{JohnsonFinite}.  We will need two facts
about $S$.  First, $S$ consists of {\em bounding pair maps}, that is, mapping classes $T_{\gamma_1} T_{\gamma_2}^{-1}$ where the $\gamma_i$ are
disjoint nonseparating curves so that $\gamma_1 \cup \gamma_2$ separates $\Sigma_g$.  Second, for
$T_{\gamma_1} T_{\gamma_2}^{-1} \in S$, either $\gamma_1 \cap v = \gamma_2 \cap v = \emptyset$ or a regular neighborhood
of $\gamma_1 \cup \gamma_2 \cup v$ is homeomorphic to the curves pictured on the left hand side of Figure \ref{figure:boundingtwist}.b.  These
facts imply that for $s \in S^{\pm 1}$, either $s \cdot v = v$ or (as demonstrated by Figure \ref{figure:boundingtwist}.b) $s \cdot v$ is
disjoint from $v$.  Condition 2 follows. \qed
\end{likeproof}

\InsertFigure{figure:boundingtwist}{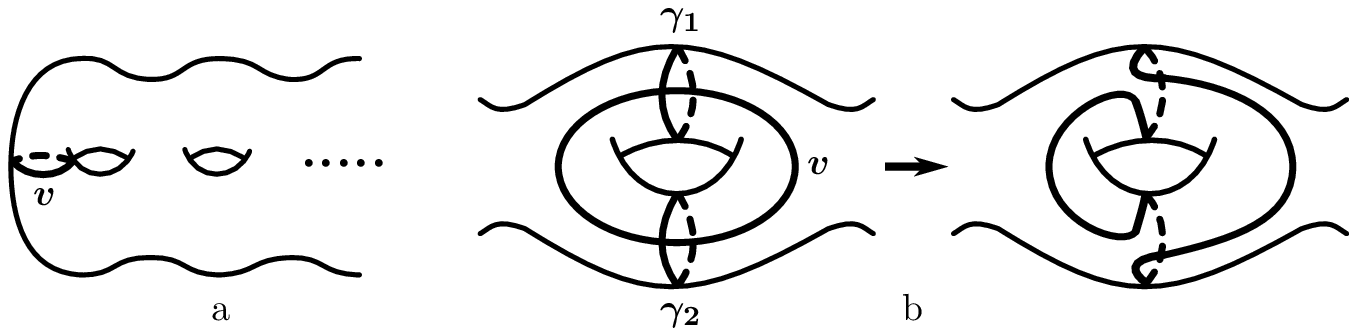}{a. Our base vertex in $\Curves^{\gamma}(\Sigma_g)$ \CaptionSpace b. $T_{\gamma_1} T_{\gamma_2}^{-1}(v)$ is disjoint from $v$}

\section{Simple connectivity}
We conclude this paper by proving Theorem \ref{theorem:csepsimplyconnected}.

\begin{likeproof}[Proof of Theorem \ref{theorem:csepsimplyconnected}]
Let 
$$S=\{T_{\alpha_1},T_{\delta_1},\ldots,T_{\alpha_g},T_{\delta_g},T_{\beta_1},\ldots,T_{\beta_{g-1}},h\}$$ 
be the collection of twists
about the curves in Figure \ref{figure:modconnected}.a together with the hyperelliptic involution $h$ (see \cite[page 52]{IvanovSurvey} for
the definition of $h$; the need for $h$ will become clear shortly).  Also, let $v$ and
$v'$ be the curves in Figure \ref{figure:modconnected}.b.  There is a natural map $\Mod(\Sigma_g) \rightarrow \CSep(\Sigma_g)$ taking $g$ to
$g(v)$.  Closely examining the proofs of Lemma \ref{lemma:main} and Theorem \ref{theorem:curvesconnected}, we see that they say that this map 
extends to a $\Mod(\Sigma_g)$-equivariant map
$$\phi : \Cayley(\Mod(\Sigma_g),S) \longrightarrow \CSep(\Sigma_g).$$
Here $\Cayley(\Mod(\Sigma_g),S)$ is the {\em Cayley graph} of $\Mod(\Sigma_g)$, that is, the graph whose vertices are elements of
$\Mod(\Sigma_g)$ and where $g_2$ is connected by an edge to $g_1$ if $g_2 = g_1 s$ for some $s \in S$.  We will prove that
the induced map $\phi_{\ast} : \pi_1(\Cayley(\Mod(\Sigma_g),S),1) \rightarrow \pi_1(\CSep(\Sigma_g),v)$ is the zero map by examining
the images of the loops associated to a set of relations for $\Mod(\Sigma_g)$.  We will then show that $\phi_{\ast}$ is surjective; this will allow
us to conclude that $\pi_1(\CSep(\Sigma_g),v)=0$, as desired.

\begin{liketheo}[Claim 1]
$\phi_{\ast} : \pi_1(\Cayley(\Mod(\Sigma_g),S),1) \rightarrow \pi_1(\CSep(\Sigma_g),v)$ is the zero map.
\end{liketheo}

\begin{likeproof}[Proof of claim]
It is well known that we can construct a simply connected complex $X$ from $\Cayley(\Mod(\Sigma_g),S)$ by attaching discs to the $\Mod(\Sigma_g)$-orbits
of the loops associated to any complete set of relations for $\Mod(\Sigma_g)$.  We will show that
the images in $\CSep(\Sigma_g)$ of the loops associated to these relations are contractible.  This will imply that we can
extend $\phi$ to $X$.  Since $X$ is simply connected, we will be able to conclude that $\phi_{\ast}$ is the zero map, as desired.

Now, the loop in $\Cayley(\Mod(\Sigma_g),S)$
associated to a relation $s_1 \cdots s_k = 1$ with $s_i \in S^{\pm 1}$ is 
$1 - s_1 - s_1 s_2 - \ldots - s_1 s_2 \cdots s_k$.  Since the only elements of $S^{\pm 1}$ which act
non-trivially on $v$ are $T_{\beta_{g-1}}^{\pm 1}$, the function 
$\phi$ maps the edge $s_1 \cdots s_{i-1} - s_1 \cdots s_i$ to a fixed vertex unless $s_i = T_{\beta_{g-1}}^{\pm 1}$,
in which case it maps it to the path $s_1 \cdots s_{i-1} (v) - s_1 \cdots s_{i-1} (v') - s_1 \cdots s_i(v)$.  Observe that
the only elements of $S^{\pm 1}$ which act non-trivially on $v'$ are $T_{\beta_1}^{\pm 1}$.  If
none of the $s_i$ equal $T_{\beta_1}^{\pm 1}$, then for all $i$ we would have $s_1 \cdots s_{i-1} (v')=v'$, so $\phi$ would
take the loop associated to the relation $s_1 \cdots s_k=1$ to a loop of the form $v_1 - v_2 - \ldots - v_{2l+1}$ with $v_{2i} = v'$ for $1 \leq i \leq l$.  
This loop can be contracted to $v'$.

We therefore only need to worry about relations which involve both $T_{\beta_1}^{\pm 1}$ and $T_{\beta_{g-1}}^{\pm 1}$.  By
Theorem \ref{theorem:modpresentation} from the appendix, we can find a presentation for $\Mod(\Sigma_g)$ whose
generators are $S$ and whose only relations involving both $T_{\beta_1}^{\pm 1}$ and $T_{\beta_{g-1}}^{\pm 1}$ are
\begin{equation}
T_{\beta_1}^{-1} T_{\beta_{g-1}}^{-1} T_{\beta_1} T_{\beta_{g-1}} = 1 \label{commute}
\end{equation}
and
\begin{equation}
T_{\alpha_g} T_{\beta_{g-1}} T_{\alpha_{g-1}} \cdots T_{\beta_1} T_{\alpha_1} T_{\delta_1}^2 T_{\alpha_1} T_{\beta_1} \cdots T_{\beta_{g-1}} T_{\alpha_g} h^{-1} = 1. \label{hyper}
\end{equation}
We conclude that we must only check that the $\phi$-image of the loops associated to these two relations are contractible.  For
the relation given in (\ref{commute}), it is clear that we can find a separating curve disjoint from every vertex of the associated
$\CSep(\Sigma_g)$-loop (for instance, $v''$ in Figure \ref{figure:modconnected}.b is such a curve), so this loop is contractible.  For
the relation given in (\ref{hyper}), the associated loop is the following, where we suppress the (trivial) edges $w(v) - w s (v)$ associated
to generators $s \in S^{\pm 1}$ not equal to $T_{\beta_{g-1}}^{\pm 1}$:
\begin{align*}
&v - T_{\alpha_g}(v') - T_{\alpha_g} T_{\beta_{g-1}}(v) \\
&\ \ \ \ \ \ - T_{\alpha_g} T_{\beta_{g-1}} T_{\alpha_{g-1}} \cdots T_{\beta_1} T_{\alpha_1} T_{\delta_1}^2 T_{\alpha_1} T_{\beta_1} \cdots T_{\alpha_{g-1}}(v') - v \\
= &v - v' - T_{\alpha_g} T_{\beta_{g-1}}(v) - v' - v.
\end{align*}
This is clearly contractible, so the claim follows. \qed
\end{likeproof}

\begin{liketheo}[Claim 2]
$\phi_{\ast} : \pi_1(\Cayley(\Mod(\Sigma_g),S),1) \rightarrow \pi_1(\CSep(\Sigma_g),v)$ is surjective.
\end{liketheo}

\begin{likeproof}[Proof of claim]
We first find a sufficient condition for a loop to lie in the image of $\phi_{\ast}$.  Consider any
loop $\ell = v_1 - v_2 - \ldots - v_{2n+1}$ in $\CSep(\Sigma_g)$ with $v_1 = v_{2n+1}=v$.  Assume that
each $v_i$ is a separating curve which cuts off a 1-holed torus and that for $0 \leq i <n$ there
exists a simple closed nonseparating curve $\epsilon_i$ and some $e_i = \pm 1$ so that 
$v_{2i+3} = T_{\epsilon_i}^{e_i} (v_{2i+1})$, so that $i_g(\epsilon_i, v_{2i+1}) = 2$, and
so that $i_g(\epsilon_i, v_{2i+2}) = 0$ (for instance, we could have $v_{2i+1}=v$, $v_{2i+2}=v'$, and $\epsilon_i = \beta_{g-1}$).  
We claim that $\ell$ is in the image of $\phi_{\ast}$.  

To begin with, it is enough to find some word $w$ in $S^{\pm 1}$ (not necessarily a
relation) so that $\ell$ is the image under $\phi_{\ast}$ of the path in $\Cayley(\Mod(\Sigma_{g}),S)$ associated to $w$.  Indeed,
we then would have $w(v) = v$.  Since $S \setminus \{T_{\beta_{g-1}}\}$ generates the stabilizer in $\Mod(\Sigma_g)$
of $v$, we can find some word $w'$ in $(S \setminus \{T_{\beta_{g-1}}\})^{\pm 1}$ so that $w w' =1$; this is the desired
relation.

We will prove the existence of $w$ by induction on $n$ (in this part of the proof, we do not assume that $\ell$ is a loop).  
The case $n=0$ being trivial, we assume that $n>0$.  Using the induction
hypothesis, we can find a word $w_{n-1}$ so that $\phi_{\ast}$ takes the path associated to $w_{n-1}$ to $v_1 - \ldots - v_{2n-1}$.  Observe 
that
\begin{align*}
&i_g(w_{n-1}^{-1}(v_{2n}),w_{n-1}^{-1}(\epsilon_n)) = i_g(v_{2n},\epsilon_n) = 0, \\
&i_g(w_{n-1}^{-1}(v_{2n}), v) = i_g(v_{2n},w_{n-1}(v)) = i_g(v_{2n},v_{2n-1}) = 0, \\
&i_g(w_{n-1}^{-1}(\epsilon_n), v) = i_g(\epsilon_n,w_{n-1}(v)) = i_g(\epsilon_n, v_{2n-1}) = 2.
\end{align*}
This implies that there must exist some $f \in \Mod(\Sigma_g)$ so that $f(v') = w_{n-1}^{-1}(v_{2n})$, 
so that $f(\beta_{g-1}) = w_{n-1}^{-1}(\epsilon_n)$, and so that $f(v) = v$.  Since the stabilizer in $\Mod(\Sigma_g)$
of $v$ is generated by $S \setminus \{T_{\beta_{g-1}}\}$, we can find some word $w''$ in $(S \setminus \{T_{\beta_{g-1}}\})^{\pm 1}$
so that $w'' = f$.  We claim that $w=w_{n-1} w'' T_{\beta_{g-1}}^{e_n}$ works.  Indeed, since $w_{n-1} f(v') = v_{2n}$, 
the path associated to $w$ consists of $v_1 - \ldots - v_{2n-1}$ followed by the path
\begin{align*}
w_{n-1} f (v') - w_{n-1} f T_{\beta_{g-1}}^{e_n} (v) &= v_{2n} - w_{n-1} (f T_{\beta_{g-1}}^{e_n} f^{-1}) (v) \\
&= v_{2n} - w_{n-1} T_{f (\beta_{g-1})}^{e_n} (v) \\
&= v_{2n} - w_{n-1} T_{w_{n-1}^{-1} (\epsilon_n)}^{e_n} (v) \\
&= v_{2n} - w_{n-1} w_{n-1}^{-1} T_{\epsilon_n}^{e_n} w_{n-1} (v) \\
&= v_{2n} - T_{\epsilon_n}^{e_n} (v_{2n-1}).
\end{align*}
Since $T_{\epsilon_n}^{e_n} (v_{2n-1}) = v_{2n+1}$, this is the desired path.

Now consider an arbitrary $\ell' \in \pi_1(\CSep(\Sigma_g),v)$.  We claim that we can homotope $\ell'$ so that it satisfies the
above condition.  In fact, we will prove more generally that if $\ell'$ is any (not necessarily closed) path starting at $v$ whose
final endpoint corresponds to a curve cutting off a 1-holed torus, then we
can homotope it (fixing the endpoints) so that it satisfies all of the above conditions except for the closedness of
the path.  

We can assume without loss of generality that $\ell'$ is a simplicial path in the 1-skeleton.  It is an
easy exercise to see that we can homotope $\ell'$ so that all of its vertices cut off 1-holed tori, and
in addition we can arrange for $\ell'$ to contain an odd number of vertices and for no two adjacent vertices of $\ell'$ to be
identical.  Enumerate the vertices of 
of $\ell'$ as $v_1 - \ldots - v_{2m+1}$.  By induction on $m$, we can assume that $v_1 - \ldots - v_{2m-1}$ satisfies
the desired condition.  Now, using standard properties of $\Mod(\Sigma_g)$ we can find a sequence of
simple closed curves $\eta_1,\ldots,\eta_k$ and numbers $f_1,\ldots,f_k=\pm 1$ so that
$T_{\eta_1}^{f_1} \cdots T_{\eta_k}^{f_k} (v_{2m-1}) = v_{2m+1}$ and so that for $1 \leq i \leq k$ we have
$i_g(\eta_i,v_{2m-1})=2$ and $i_g(\eta_i,v_{2m})=0$.  We can then homotope $\ell'$ (adding ``whiskers'') so that
the path $v_{2m-1} - v_{2m} - v_{2m+1}$ is replaced by
\begin{align*}
v_{2m-1} &- v_{2m} - T_{\eta_1}^{f_1} (v_{2m-1}) - v_{2m} - T_{\eta_1}^{f_1} T_{\eta_2}^{f_2} (v_{2m-1}) \\
         &- v_{2m} - \ldots - T_{\eta_1}^{f_1} \cdots T_{\eta_k}^{f_k} (v_{2m-1}) = v_{2m+1},
\end{align*}
thus proving the claim. \qed
\end{likeproof}

\noindent
This completes the proof of Theorem \ref{theorem:csepsimplyconnected}. \qed
\end{likeproof}

\appendix
\section{Appendix : A variant on the Wajnryb presentation}

\begin{theorem}
\label{theorem:modpresentation}
For $g \geq 4$, the group $\Mod(\Sigma_g)$ has a presentation $\langle S|R \rangle$ satisfying the following conditions.
\begin{itemize}
\item $S$ is the set of Dehn twists $\{T_{\alpha_1},T_{\gamma_1},\ldots,T_{\alpha_g},T_{\gamma_g},T_{\beta_1},\ldots,T_{\beta_{g-1}}\}$ depicted in
Figure \ref{figure:modconnected}.a together with the hyperelliptic involution $h$.
\item The only relations $r \in R$ which involve both
$T_{\beta_1}^{\pm 1}$ and $T_{\beta_{g-1}}^{\pm 1}$ are
$$T_{\beta_1}^{-1} T_{\beta_{g-1}}^{-1} T_{\beta_1} T_{\beta_{g-1}} = 1$$
and
$$T_{\alpha_g} T_{\beta_{g-1}} T_{\alpha_{g-1}} \cdots T_{\beta_1} T_{\alpha_1} T_{\delta_1}^2 T_{\alpha_1} T_{\beta_1} \cdots T_{\beta_{g-1}} T_{\alpha_g} h^{-1} = 1.$$
\end{itemize}
\end{theorem}
\begin{proof}
The presentation described in this theorem is a variant of the standard Wajnryb presentation (\cite{WajnrybPresentation}; see \cite{BirmanWajnryb}
for errata).  The generating set for the Wajnryb presentation is 
$$S' = \{T_{\alpha_1}, \ldots, T_{\alpha_g}, T_{\beta_1},\ldots, T_{\beta_{g-1}}, T_{\gamma_1}, T_{\gamma_2}\}.$$
There are four families of relations.  In the notation of \cite{WajnrybPresentation}, the first three are labeled A, B, and C.  The relation
$T_{\beta_1}^{-1} T_{\beta_{g-1}}^{-1} T_{\beta_1} T_{\beta_{g-1}} = 1$ is the only relation from family A (the ``braid relations'') involving
both $T_{\beta_1}^{\pm 1}$ and $T_{\beta_{g-1}}^{\pm 1}$.  Families B (the ``two-holed torus relation'') and C (the ``lantern relation'') do
not involve both $T_{\beta_1}^{\pm 1}$ and $T_{\beta_{g-1}}^{\pm 1}$.  The final relation D (as corrected by \cite{BirmanWajnryb}) says that
the hyperelliptic involution $h$ commutes with $T_{\delta_g}$; both $h$ and $T_{\delta_g}$ are expressed using rather complicated 
formulas involving the generators $S'$.  Our relation
$$T_{\alpha_g} T_{\beta_{g-1}} T_{\alpha_{g-1}} \cdots T_{\beta_1} T_{\alpha_1} T_{\delta_1}^2 T_{\alpha_1} T_{\beta_1} \cdots T_{\beta_{g-1}} T_{\alpha_g} h^{-1} = 1$$
expands out Wajnryb's formula for $h$.  As was observed in \cite[Remark 1.a]{BirmanWajnryb}, the expression for $T_{\delta_g}$ in
terms of $S'$ used by Wajnryb is unimportant; any correct expression will work.  Now, using an argument of Humphries \cite{HumphriesGenerators},
for $1 \leq i \leq g-2$ we can express $T_{\delta_{i+2}}$ as a complicated product of elements in 
$$\{T_{\alpha_i},T_{\alpha_{i+1}},T_{\alpha_{i+2}},T_{\beta_i},T_{\beta_{i+1}},T_{\delta_i},T_{\delta_{i+1}}\}^{\pm 1}.$$ 
This allows us eliminate $T_{\delta_i}$ from $S$ for $i \geq 3$ by adding relations which do not involve both 
$T_{\beta_1}^{\pm 1}$ and $T_{\beta_{g-1}}^{\pm 1}$.
Our final relation is $[h,T_{\delta_g}]=1$; since this does not involve either $T_{\beta_1}$ or $T_{\beta_{g-1}}$, we are done. \qed
\end{proof}

\begin{likeexample}[acknowledgements]
I wish to thank my advisor Benson Farb for his enthusiasm and numerous comments, Dan Margalit
and Saul Schleimer for encouraging me to write this paper, and Matt Day and Julia Putman for offering corrections
to previous versions of this paper.  I also wish to thank that mathematics department of the Georgia
Institute of Technology for their hospitality during the time in which this paper was conceived.
\end{likeexample}

\begin{address}
Andrew Putman\\
Department of Mathematics; MIT, 2-306
77 Massachusetts Avenue
Cambridge, MA 02139-4307
\email{andyp@math.mit.edu}
\end{address}

\end{document}